\documentclass[11pt]{article}
\usepackage{setspace}
\usepackage{textcomp}
\usepackage{amsmath, amsthm, amssymb}
\usepackage[left=1.5in,top=1in,right=1in,nohead]{geometry}
\usepackage{hyperref}
\usepackage{latexsym}
\usepackage{mathrsfs}
\usepackage{rawfonts}
\usepackage[dvips]{graphicx}

\def\a{\alpha}

\def\bct{\begin{center}}
\def\ect{\end{center}}
\def\beg{\begin}

\def\<{\langle}
\def\>{\rangle}

\def\mbb{\mathbb}

\def\mc{\mathcal}

\def\tn{\textnormal}
\def\wt{\widetilde}

\newtheorem{thm}{Theorem}[section]

\title{A class of locally conformally flat 4-manifolds}

\author{Selman Akbulut\footnote{The first named author is partially supported by NSF grant DMS 9971440}
 \and Mustafa Kalafat
 }

\begin{document}
\maketitle

\begin{abstract} We construct infinite families of non-simply connected locally conformally flat (LCF) 4-manifolds realizing rich topological types. These manifolds have strictly negative
scalar curvature and the underlying topological 4-manifolds do not
admit any Einstein metrics. Such 4-manifolds are of particular interest as  examples of Bach-flat but non-Einstein spaces in the non-simply connected case. 
Besides that the underlying smooth manifolds are 
examples of spaces that admit open book decomposition in dimension 4.

\end{abstract}




\maketitle

\section{Introduction}

A Riemannian n-manifold $(M,g)$ is called {\em locally conformally flat (LCF)} if there is a function $f:U\to \mbb R^+$ in a
neighborhood of each point $p\in M$ such that $\wt{g}=fg$ is a
flat metric on $U$. It turns out that there is a simple tensorial description of this elaborate condition.
The Weyl curvature tensor 
is defined as
$$W_{ijkl}=R_{ijkl}+{R\over (n-1)(n-2)}\left|  \begin{array}{@{}cc@{}}
    g_{ik} & g_{il} \\
    g_{jk} & g_{jl} \\ \end{array} \right|
-{ 1 \over n-2 }( \left|  \begin{array}{@{}cc@{}}
    R_{ik} & g_{il} \\
    R_{jk} & g_{jl} \\ \end{array} \right|
+            \left|  \begin{array}{@{}cc@{}}
    g_{ik} & R_{il} \\
    g_{jk} & R_{jl} \\ \end{array} \right| ).$$
It is a nice exercise in tensor analysis \cite{jeffdg} that for $n\geq 4$, $M$ is LCF if and only if $W=0$. In dimension $3$ this role is taken over by the {\em Cotton} tensor, and in dimension $2$ all manifolds are LCF. The Weyl curvature tensor yields a symmetric operator $\mathcal W : \Lambda^2\to\Lambda^2$ defined by the formula $$\mathcal W (\omega)={1\over 4}W_{ijkl}\,\omega_{kl}\,e^i\wedge e^j$$
where $\{ e^1 ,.., e^n \}$ is an orthonormal basis of the $1$-forms. We are mainly concerned with dimension 4, and in this case the space of the 2-forms decomposes into the $\pm1$ eigenspaces of the Hodge star operator  $\Lambda^2=\Lambda^2_+\oplus\Lambda^2_-$.
Furthermore the operator $\mathcal W$ sends (anti-) self-dual 2-forms to
(anti-) self-dual 2-forms, hence inducing the decomposition $\mc W = \mc W^+\oplus \mc W^-$. 
 We call a Riemannian manifold $M$  {\em self-dual (SD)} if $\mc W^-=0$, and
 {\em anti-self-dual (ASD)} if $\mc W^+=0$. In these terms $M$ is LCF if and only if it is SD and ASD at the same time. 
For basics of LCF manifolds we refer to \cite{matsumoto,jeffdg}. Some common examples in dimension four are the manifolds with constant sectional curvature, 
product of two constant sectional curvature metrics of curvature
$1$ and $-1$, e.g. $S^2\times\Sigma_g ~\tn{for}~ g\geq 2$, product of a manifold of constant sectional curvature with $S^1$ or $\mbb R$. See \cite{kalafatstructuressdlcf} for a recent survey of LCF and self-dual structures on basic $4$-manifolds. 
Our main result is the following.
\begin{thm}\label{main} There are 
infinite families of closed, non-simply connected,
locally conformally flat 4-manifolds, called {\em panelled web 4-manifolds}, with Betti number growth:  $b_1\to\infty$, $b_2\to\infty$ or bounded, and $\chi\to-\infty$.
These manifolds have strictly negative scalar curvature.
\end{thm}
\vspace{.05in}
\noindent We show that many new topological types can be realized.
The idea is to conformally compactify  $S^1 \times M^3$ where $M$ is a 
hyperbolic 3-manifold with boundary. The reader will see that the resulting manifold is closed but it is 
not simply $S^1$ cross a 3-manifold. It is obtained through spinning around the
boundary of the 3-manifold. 
Recall:





\begin{thm}[\cite{braam}]\label{braamthm} Let $\bar{M}^3$ be an oriented, geometrically
finite complete hyperbolic manifold with nonempty boundary, such that $\partial \bar{M}=\cup S_j$ consists of either a disjoint union surfaces of genus $\geq 2$, or $\bar{M}=D^2\times S^1$. Let $M$ be the interior of $\tilde{M}$. Then $M\times S^1$ has a oriented closed, smooth conformal compactification $X^4$, with an $S^1$ action.

\vspace{.05in}

$X$ is locally conformally flat (LCF). 
The action has the fixed point sets conformal to the boundary surfaces $\cup S_j$ of $\bar M$ (the ideal points of the compactification). The normal bundles of the fixed
surfaces are trivial  with $S^1$ weight $1$. The hyperbolic
structure on $M$ can be recovered from $X$ by giving
$X-\cup S_j$ the metric in the conformal class for which the $S^1$ orbits have length $2\pi$. Then $M$ is the Riemannian quotient of $X-\cup S_j$ by $S^1$. \end{thm}
\noindent In particular the connected sums  $\sharp_n S^3\times S^1$ and $S^2\times \Sigma_g$ for $g\geq 2$
can be obtained from this theorem. 
In the first case we begin with several cyclic groups of isometries of $\mbb H^3$ each of which
yields a quotient $D^2\times S^1$, combining them by the first combination theorem gives a classical Schootky group corresponding to the boundary connected sums
of the corresponding $D^2\times S^1$s. Boundary connected sum in three dimensions corresponds to the ($S^1$ equivariant conformal)
connected sum in four dimensions.
In the second case we begin with a Fuchsian group of isometries of
$\mbb H^3$, yields a quotient $I\times \Sigma_g$.

\vspace{.05in}

In this paper we begin with a more general class of Kleinian groups
called the {\em panelled web groups}, constructed by Bernard Maskit in \cite{maskitpwg}. After the application of the Theorem \ref{braamthm}, we obtain 4-manifolds with more complicated topology. 
We describe concrete handlebody pictures of these manifolds in terms of framed links, which describes their smooth topology. We will call these LCF manifolds {\em panelled web 4-manifolds}\index{LCF 4-manifolds, Panelled Web}. We hope that our concrete ``visual'' techniques here will be useful in constructing special metrics on other manifolds, especially the other no-simply connected ones. 

\vspace{.05in}

We are also able to compute the sign of the scalar curvature for the panelled web 4-manifolds. Recall that by the solution of the Yamabe problem, any Riemannian metric on a closed manifold is conformally equivalent to the one with constant scalar curvature. And the sign of this constant is an invariant of the conformal structure, called the {\em type} of the metric or its conformal class. Using the results of \cite{cltopsd}  and additionally \cite{schoenyau88,nayatani97} we can show the following. 

\beg{thm}\label{signscalarcurvaturenegative}
The conformal class of the natural metric on the panelled web 4-manifolds is of negative type, i.e. the metric can be rescaled to have constant negative scalar curvature. In the case of $b_2\neq 0$, more generally the underlying topological manifolds of panelled web 4-manifolds do not admit any locally conformally flat metric of positive or zero scalar curvature. 
\end{thm}
\noindent 
Considering the natural metric of these manifolds, one can also directly compute its sign through the Hausdorff dimension of the Kleinian groups 
used to uniformize the related hyperbolic 3-manifold. See the Section 
\S\ref{secsign} for the details. 

\vspace{.05in}

Finally, we can give an answer to the problem of whether the underlying  smooth 4-manifolds admit any Einstein metric. 
We compute the Euler characteristics of the manifolds we construct. The Euler characteristics of the building blocks are all strictly negative, since the Euler characteristic is additive, and it turns out to
be strictly negative for all of our panelled web 4-manifolds. By  the
generalized Gauss-Bonnet Theorem we express the Euler characteristic
$\chi$ of a 4-manifold as
$$\chi (M) = \frac{1}{8\pi^2} \int_M ~~\frac{s^2}{24}-\frac{|\stackrel{\circ}{\tn{Ric}}|^2}{2}+
 |W|^2~~dV_g.$$
If $M$ admits an Einstein metric, then the trace free Ricci
curvature tensor $$\stackrel{\circ}{\tn{Ric}}\hspace{.1cm}=\tn{Ric}-{s\over
4}g$$ vanishes identically. So that $\chi\geq0$, which implies the
following

\beg{thm} The topological manifolds underlying the panelled
 web 4-manifolds do not admit any Einstein metrics.\end{thm}
\noindent 
This is interesting because of the following.  Einstein metrics are have vanishing Bach tensor, so that they are Bach-flat (BF). LCF metrics are also BF. Then our examples are BF but not Einstein. Therefore, in the 
highly non-simply connected case, these examples illustrates the converse statement.  
See also 6.32 of \cite{besse} for simpler examples.  
It is easier to give simply-connected examples of this phenomenon; 
$\sharp_n\mbb{CP}_2$ carries self-dual metrics by $\cite{clexplicit}$ however no Einstein metric for $n\geq 4$ by the Hitchin-Thorpe inequality.

\vspace{.03in}

It is a curious question whether these smooth manifolds
carry any optimal metric \cite{cloptimal}. 
Since they do not admit any Einstein metric, the first possibility is eliminated. Another possibility of being scalar-flat anti-self-dual (SF-ASD) can also be eliminated in $b_2\neq 0$ case, since the techniques mentioned in Section 
\S\ref{secsign} goes through in this case as well. Besides that, since the signature of these manifolds vanish, self-duality or anti-self-duality of the metric is equivalent to being locally conformally flat in this case. 
Consequently the optimal metric problem currently remains open for these manifolds.  




\vspace{.05in}

Note that the handlebody pictures are 
essential to deal with non-simply connected
manifolds in general. This is the standard and only way to define and understand them generally. 
Otherwise one trapped into products and connected sums.
There is no way to get complicated topological types other than
showing the explicit surgery scheme. They are somehow the definition of the manifolds.
Products of simple manifolds and their connected sums constitute a set of measure zero in the whole family of non-simply connected $4$-manifolds. 
Because of this reason, 
we consider this study as a foundational work to analyze, give
examples of LCF (and also SD) metrics on non-simply
connected spaces. 
This work has many further applications. In a forthcoming paper \cite{ak1} using the techniques here, we construct self-dual but not locally conformally flat metrics on families of non-simply connected 4-manifolds with small signature. Secondly, in \cite{ak2} we analyze the existence of symplectic, almost complex and complex structures on the panelled web 4-manifolds constructed here, and give interesting counterexamples. 
More applications are on the way.

\vspace{.05in}

In section \S\ref{secpanelled} we review the
hyperbolic 3-manifolds which we use in our constructions. In \S\ref{sechandlebody} we
describe the topology of the building blocks of the 4-manifolds in
interest, by constructing their handlebody pictures. In \S\ref{secinvariants} we compute the algebraic topological invariants of these 4-manifolds. In \S\ref{secsequences} we construct interesting sequences of locally conformally flat 4-manifolds
by using these building blocks. Finally in
\S\ref{secsign} we compute the sign of the scalar curvature of the
metrics on these manifolds.

\vspace{.05in}

{\bf Acknowledgments.} The second author would like to thank to
Claude LeBrun for his suggestions, Bernard Maskit, for generously sharing his knowledge, and many thanks to Alphan Es, Jeff Viaclovsky, Yat-Sun
Poon, Caner Koca, Feza G\"ursey Institute members at \.Istanbul. The figures are
sketched by the program IPE of Otfried Cheong. We thank Selahi
Durusoy for helping with the graphics. Also thanks to the anonymous  referee for many useful comments.

\section{Panelled Web Groups}\label{secpanelled}

In this section we will describe the 3-manifolds from which we
construct our LCF 4-manifolds. These are closed hyperbolic 3-manifolds, which are obtained by
dividing out the hyperbolic 3-space $\mathbb H^3$ with a group of
its isometries. The isometry group is a discrete group obtained out
of certain Fuchsian and extended-Fuchsian groups, by taking their
combinations using the theorems of Maskit.
In 1981 B. Maskit introduced this new class of Kleinian groups called the {\em panelled web groups}, and gave a set of examples. Here we  first review the constructions in
\cite{maskitpwg}.

\beg{defn}
 A {\em Fuchsian group} is a discrete
group of fractional linear transformations $z\mapsto (az+b)/(cz+d)$
acting on the hyperbolic plane\footnote{We will be using the upper
half plane model of the hyperbolic plane throughout this paper.}
$\mathbb H^2$, where $ad-bc\neq 0$ and $a,b,c,d$ are real. The group
is of the {\em first kind} if every real point is a limit point, it
is of the {\em second kind} otherwise.
\end{defn}

M\"obius transformations can be written as a composition of
reflections and inversions. These motions act on the extended
complex line $\mathbb{\hat{C}}$ as well as on the upper half space
$\mathbb H^3=\{(z,t)|z\in\mathbb C , t\in\mathbb R^+\}$ by the
usual way.
In our case the trasformations preserves the $\mathbb H^2$ so that
they are written as a product of reflections and inversions in lines
and circles which are orthogonal to the real line. The extended
motions in $\mathbb H^3$ preserve the planes passing through the
real line, it follows that if $G$ is a Fuchsian group then, $\mbb
H^3/G=\mbb H^2/G\times (0,1)$.

\vspace{.05in}

A group of M\"obius transformations is called {\em
elementary}\index{elementary} if it has at most two limit points. As
an example, a hyperbolic cyclic group $H=\<z\mapsto \lambda^2 z\> ,
\lambda\neq 1$ or its conjugates has two limit points and $\mbb
H^2/H$ is an annulus. Another is a trivial group, it has no limit
point and $\mbb H^2/\{1\}$ is a disk.
Let $\Sigma_{g,n}$ be the interior of a compact orientable surface
with boundary, where $g$ and $n$ stand for the genus and number of
boundary components, respectively. Assume $\Sigma_{g,n}$ is neither a disk nor an
annulus. Then there is a purely hyperbolic, non-elementary Fuchsian
group of the second kind $G$ so that $\mbb H^3/G =
\Sigma_{g,n}\times (0,1)$. Conversely, if $G$ is a finitely
generated, purely hyperbolic, non-elementary Fuchsian group of the
second kind, then $\mbb H^2/G$ is the interior of a compact
orientable surface with boundary neither a disk nor an annulus, so
that $\mbb H^3/G = \Sigma_{g,n}\times (0,1)$.

\vspace{.05in}

We can construct the group $G$ corresponding to the surface of genus
$g$ with $n$ boundary components using $4g+2(n-1)$ disjoint,
identical circles $C_1,C_1' \cdots $ $C_{2g+n-1},$ $C_{2g+n-1}'$
centered at the real line. The generators of $G$ will be M\"obius
transformations $a_i$ mapping $C_i$ to $C_i'$, which can be
constructed as a composition of an inversion in $C_i$ followed by a reflection in the perpendicular bisector of the centers of the two
circles. Using either of the combination theorems, we see that the group $G$
generated by $a_1\cdots\a_{2g+n-1}$  is discrete, and acts
freely on $\mbb H^2$.
\begin{figure}[!h] \bct
\includegraphics[width=\textwidth]{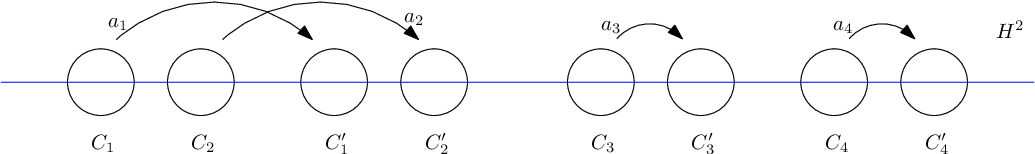}\\ \ect
  \caption{{\em {\small Schottky generators for the Fuchsian group for $\Sigma_{1,3}$.}}}
  \label{schottkygenerators}
\end{figure}
Figure \ref{schottkygenerators} shows the case for $g=1$ and $n=3$.
Notice that each generator $a_3, a_4$ generates a hole, on the other
hand the generators producing the genus $a_1 , a_2$ altogether
generates only one hole as they stick all the nearby boundary
components together.
The quotient $\mathbb H^3/G$ is the product $\Sigma \times (0,1)$ is the interior of $\Sigma \times I$ for $I=[0,1]$ which is called an
{\em I-bundle of type (i)} or a {\em trivial I-bundle} on $\Sigma$.
If there is an orientation {\em reversing}, free, involutive
homeomorphism $h:\Sigma\to\Sigma$, we extend $h$ to an orientation
{\em preserving} homeomorphism $$h':\Sigma\times I\to\Sigma\times I
 \hspace{.5cm} \tn{by} \hspace{.5cm} h'(x,t)=(h(x),1-t),$$ then we call the quotient
$\Sigma\times I/h'$ to be an {\em I-bundle of type (ii)} or a {\em
twisted I-bundle associated to $\Sigma$} or {\em over $\Sigma/h$}.
Next we will construct the Kleinian groups corresponding to the
twisted I-bundles.

\beg{defn}
A non-elementary Kleinian group which
is not itself Fuchsian, but contains a subgroup of index 2 which is
Fuchsian, is called an {\em extended Fuchsian group}\index{extended
Fuchsian group}. \end{defn}

A M\"obius transformation is called {\em parabolic, loxodromic} or
{\em elliptic} if the number of its fixed points in
$\overline{\mathbb H}^3$ is one, two or infinity, respectively. {\em
Hyperbolic} elements are the transformations conjugate to $z\mapsto
\lambda z , \lambda > 1$, which are also loxodromic. Besides, a
transformation is {\em elliptic} iff it has a fixed point in
$\mathbb H^3$.

\vspace{.05in}

If we start with a finitely generated, non-elementary, purely
loxodromic extended Fuchsian group $G$, we can write
 $G=\<g,G^0\>$, for some
Fuchsian group $G^0$, so that $g\,G^0g^{-1}=G^0$ and $g^2\in G^0$ (\cite{maskitpwg,maskitkleinian,matsuzaki}).
After renormalizing we can assume that $g$ has fixed points at
$0,\infty$ and then $g$ maps a Euclidean plane passing through the
real line with an inclination of $\alpha$ with the upper half plane
onto a Euclidean plane also passing through the real line with
inclination of $\pi - \alpha$ degrees. The plane with $\alpha=\pi/2$
is kept invariant. $G$ has no elliptic elements so it is
torsion-free, 
implying that the action of $g$ on
the $\alpha=\pi/2$ plane can have fixed points only on the real
axis. We conclude that $\mathbb H^3/G$ is equal to the $\mathbb
H^3/G^0$ modulo the action of $g$, so is an I-bundle of type (ii)
over $\mathbb H^2/G$.

\vspace{.05in}

\begin{figure}[!h] \bct
\includegraphics[width=.6\textwidth]{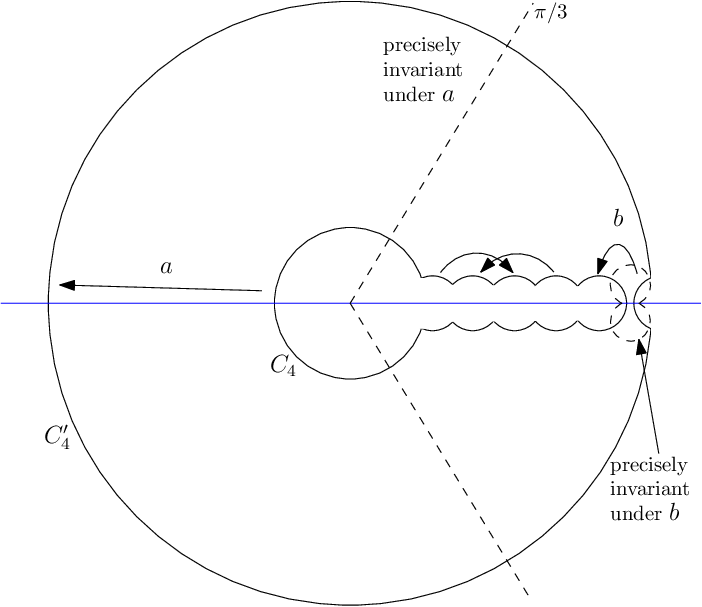}\\ \ect
  \caption{{\em {\small Fundamental region of $G_1$ as a Kleinian group.}}}
  \label{fundamentalregion1}
\end{figure}

To construct our 3-manifolds, we  glue the hyperbolic
3-manifolds obtained out of the quotients of Fuchsian and
extended-Fuchsian groups. The gluing is done along the cylinders. If
we begin with the case $n>0$, i.e. surfaces with holes, then the
quotient 3-manifolds have cylinders along the boundary,
corresponding to the boundary curves. These are of the form $W\times
I$ for a boundary curve $W$. Each boundary cylinder has a {\em
median} $W\times \{{1/2}\}$ on it, which divides it into two {\em
half cylinders}. The gluing procedure is to glue these half
cylinders by the standard homeomorphism matching the medians to get
a connected 3-manifold at the end, which does not have any more
spare (unglued) half cylinders. Then we finish the construction
with the optional complex twist operation along some of the medians.
All of these operations are done using the combination theorems,
which never lead us out of the class of geometrically finite groups. Gluing the half cylinders of two different
3-manifolds is achieved by the following:
 \beg{thm}[First
Combination \cite{combination1,combination3}] Let $G_1$ and $G_2$ be
Kleinian groups with a common subgroup $H$. Let $C$ be a simple
closed curve dividing $\hat{\mathbb C}$ into the topological disks
$B_1 , B_2$ where $B_i$ is precisely invariant under $H$ in $G_i$.
Then the group $G$ generated by $G_1$ and $G_2$ is discrete, and  $G$ is
the free product of $G_1$ and $G_2$ with amalgamated subgroup $H$.
If $D_i$'s are fundamental domains for $G_i$'s, where $D_i\cap B_i$
is a fundamental domain for the action of $H$ on $B_i$, then
$D_1\cap D_2$ is a fundamental domain for $G$.
\end{thm}
Here, a subset $A$ of $\hat{\mathbb C}$ is said to be {\em precisely invariant}
under the subgroup $H$ in $G$, if $h(A)=A$ for every $h\in H$ and $g(A)\cap A=\emptyset$ for every
$g\in G\backslash H$. 
Let us illustrate this gluing with an example from
\cite{maskitpwg} with a Fuchsian group $G_1$ and an
extended-Fuchsian group $G_2$, which will correspond to the trivial
and twisted  I-bundles  over $\Sigma_{1,2}$. Here $G_1$ is generated
by the elements whose actions are described by the circles
$C_1,C_1'\cdots C_4,C_4'$. We choose the circles generating the
genus closer to each other so that they do not generate an extra
hole, this reduces the number of boundary circles to two. We label
the elements generating these holes as $b$ and $a$, and slide the
center of the circle $C_4'$ to the right on the real axis till it
reaches $+\infty$ and then slide back from $-\infty$ to the right
till it reaches to the origin. So that the outside of $C_4$ is
mapped inside of $C_4'$ contrary to the standard mapping in Figure
\ref{schottkygenerators}.
The fundamental region of $G_1$ as a Kleinian group looks like
Figure \ref{fundamentalregion1}. $C_4'$ is the large and $C_4$ is the small circle centered at the origin.
By our choice of the circle $C_4'$ we intend to provide the common
subgroup to be $H=\<a\>$ where $a : z\mapsto \lambda z , \lambda
>1$. $a$ is a dilation which is still a schottky generator. The
dotted lines and circles denote the {\em lens angle} for $a$ and $b$,
which is the smallest angle between the real axis and the largest
precisely invariant circular region bounded above by a circle
passing through the fixed points of the group, and below by the real
axis. It is denoted by $\varphi_H$. Incidentally, $a$ and $b$ are
the {\em boundary elements} of this Fuchsian group, e.g. the generators of the hyperbolic cyclic subgroups of
a Fuchsian group of the second kind keeping invariant the segment of the real axis on which the group acts discontinously. The dashed
circles encloses invariant regions for the boundary elements $a$ and
$b$. The two lines stand for the parts of circles at infinity.

\begin{figure}[!h] \bct
\includegraphics[width=.7\textwidth]{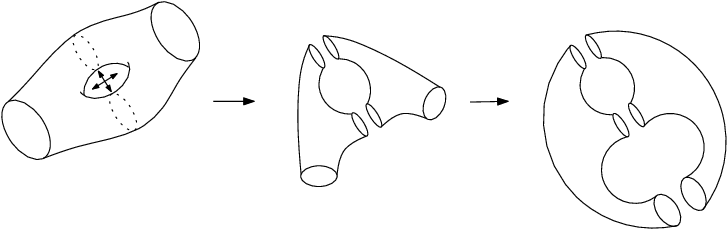}\\ \ect
  \caption{{\em {\small $\Sigma_{1,2}$ with its involution and how it sits in the fundamental region
  for $G_2$.}}}
  \label{involutionandcuttedtorusforg2}
\end{figure}

The fundamental region for $G_2$ is constructed in a more complicated way. We begin with the Fuchsian
group generating $\Sigma_{0,3}$, such that one of the holes is generated by the same $a$ as in $G_1$.
We than add a new generator $g_2$
 mapping the rest of the holes to one another.
 Adjoining this new element $g_2$ can be considered as an application of the second combination Thm \ref{secondcombination}. $G_2$ corresponds to the twisted I-bundle over $\Sigma_{1,2}$.

\begin{figure}[!h] \bct
\includegraphics[width=.6\textwidth]{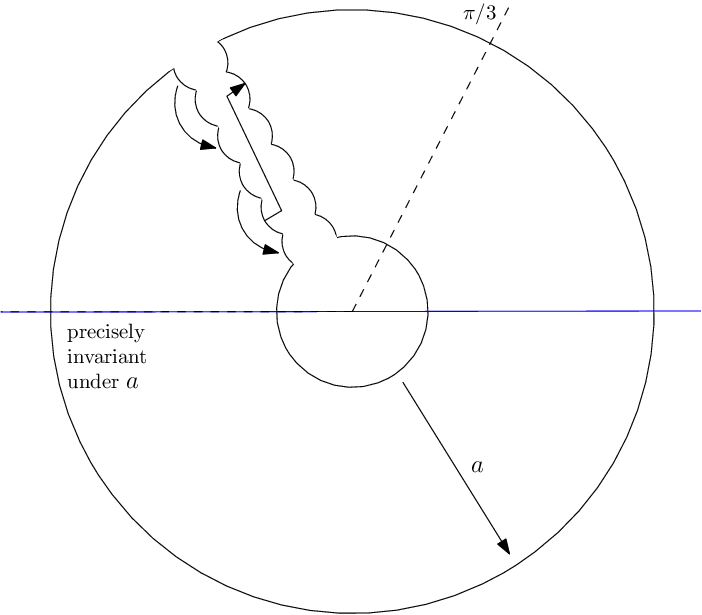}\\ \ect
  \caption{{\em {\small Fundamental region of $G_2$ as a Kleinian group.}}}
  \label{fundamentalregion2}
\end{figure}

\vspace{.05in}

\begin{figure}[!h] \bct
\includegraphics[width=.6\textwidth]{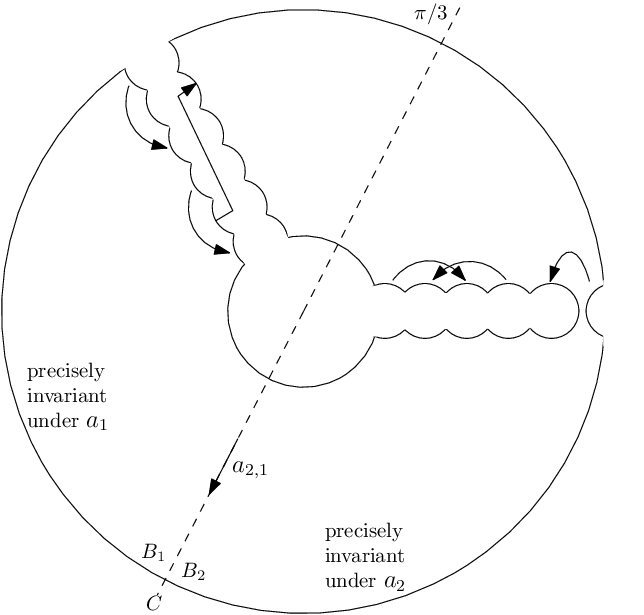}\\ \ect
  \caption{{\em {\small Fundamental region of the first combination of $G_1$ and $G_2$ along $\<a\>$.}}}
  \label{fundamentalregion3}
\end{figure}

Finally, we conjugate the group by   $g:z\mapsto exp(2\pi i/3)z$
to rotate the fundamental region by $\pi/3$ in the counter clockwise direction so that
the fixed points, geodesics of the elements of $G_2$ generated by other than $a$
lies on the other side of the line $C:\theta=\pi/3$, as in Figure \ref{fundamentalregion2} . We direct the reader to \cite{maskitpwg} for details. 
To apply the combination theorem, we take the line $C$ as the seperating circle which
seperates $\hat{\mathbb C}$ into the disks $B_1, B_2$ lying
on the left and right hand side in the Figure \ref{fundamentalregion3}, respectively. We choose our lens angles $\varphi < \pi/3$ so that $B_i$ is precisely invariant under $H=\<a\>$ in $G_i$.
The combination theorem says that the group generated by $G_1$ and $G_2$ is discrete. A fundamental domain
is as in Figure \ref{fundamentalregion3}.

In three dimension, we glued the cylinder of the twisted I-bundle to a cylinder of the
trivial I-bundle along $L/H$ where $L$ is the geodesic plane in $\mathbb H^3$ with boundary $C$.
However we only want to glue the half-cylinders. We can take apart the
glued half-cylinders and glue back in a different way using the second combination theorem.

\beg{thm}[Second Combination \cite{combination2,combination3}]\label{secondcombination} Let $G$ be a
Kleinian group with subgroups $H_1$ and $H_2$. Let $B_1 , B_2$ be two disjoint topological disks
where $(B_1,B_2)$ is precisely invariant under $(H_1,H_2)$ pairwise. Suppose there is a M\"obius
transformation $f$ mapping the interrior of $B_1$ onto the exterrior of $B_2$, where $fH_1f^{-1}=H_2$.
Then the group $G^*$ generated by $G$ and $f$ is discrete, has the relations of $G$ and $fH_1f^{-1}=H_2$.
A fundamental domain is given by $D\cap ext(B_1)\cap ext(B_2)$, where $D$ is a fundamental domain for $G$.
\end{thm}

\begin{figure}[!h] \bct
\includegraphics[width=.6\textwidth]{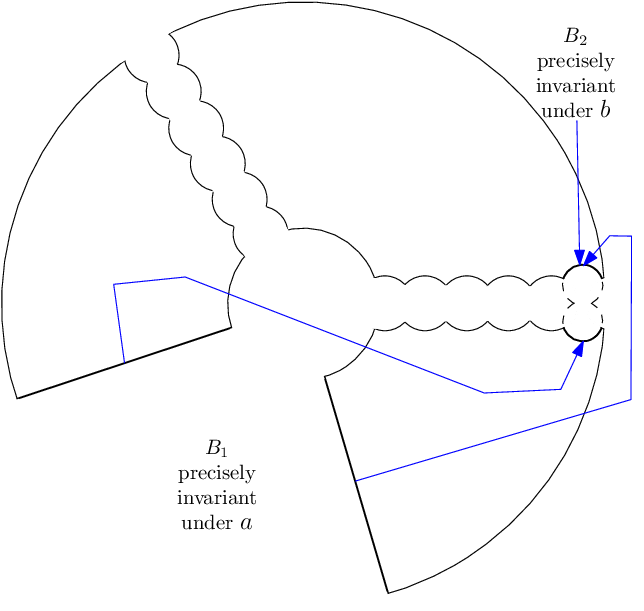}\\ \ect
  \caption{{\em {\small Application of the second combination theorem.}}}
  \label{fundamentalregion4}
\end{figure}

Here, the  {\em pairwise precise invariance} of $\{ A_1, A_2 \}$ means the usual invariance with the condition that
$gA_i\cap A_j=\emptyset$ for $i\neq j$ and for any $g\in G$.
We apply this theorem to the subgroups $\<a\>$ and $\<b\>$
in the group $G$, which we have constructed above. We arrange the loxodromic transformations
$a$ and $b$ such that they are conjugate to the transformation $z\mapsto \lambda z$
with the same $\lambda$ called the {\em multiplier}, so that they are conjugate to each other.
Choose $B_1$ as the sector $|\arg z - 4\pi/3| < \varphi$ where $\varphi<\pi/3$.
It is clearly precisely invariant under $H_1=\<a\>$ in $G$.
We choose $B_2$ to be the inside of the circular arcs passing through the
fixed point of the group $H_2=\<b\>$. We take out the sector and inside the circular arcs,
and glue the boundaries by the theorem. See Figure \ref{fundamentalregion4}.

\begin{figure}[!h] \bct
\includegraphics[width=\textwidth]{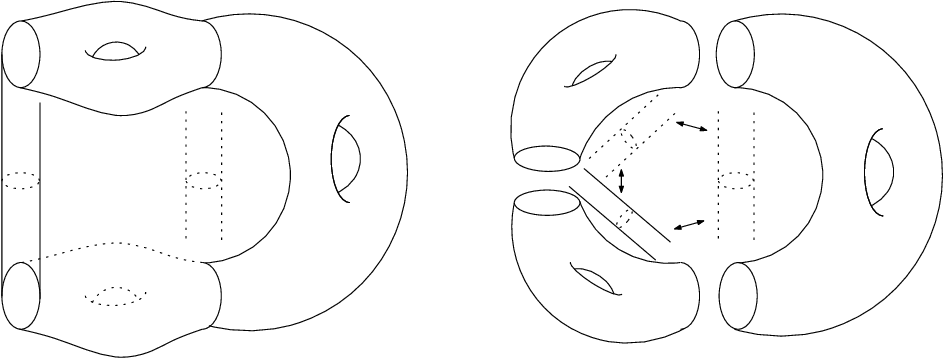}\\ \ect
  \caption{{\em {\small Effects of the first and second combination theorems in 3-dimensions.}}}
  \label{threedimension}
\end{figure}

In three dimensions, recall that applying the first combination,
we have glued a cylinder of the trivial I-bundle to the cylinder of the twisted I-bundle.
Application of the second combination tears apart one of these glued half-cylinders, and glues the
half-cylinder of the trivial I-bundle to its opposite half-cylinder, glues the spare
half-cylinder of the trivial I-bundle to the spare half-cylinder of the twisted I-bundle.
Figure \ref{threedimension} shows the identifications before and after the application of the
second combination theorem.

\vspace{.05in}

\begin{figure}[!h] \bct
\includegraphics[width=.6\textwidth]{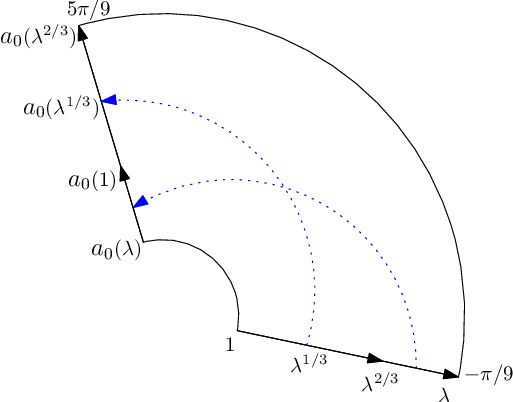}\\ \ect
  \caption{{\em {\small A fundamental region for $H_0$.}}}
  \label{annulusoverthree}
\end{figure}

Our final operation is the $p/q$ complex twist operation for relatively prime integers $p$ and $q$.
We illustrate the case for $p/q=1/3$. This will be nothing but the application of the Second Combination
Theorem to $G$ and $H_0=\<a_0\>$, where $a_0 : z\mapsto\lambda^{1/3}\exp(2\pi i/3)z$ and the common subgroup
is taken to be $H_1=\<a\>$, where $a:z\mapsto \lambda z , \lambda>1$. If we consider the isomorphism
$H_0\approx \mathbb Z$, then $H_1$ will correspond to the $3\mathbb Z$ in $\mathbb Z$ since $a_0^3=a$.
A fundamental region in $\hat{\mathbb C}$ for $H_1$ is an annulus of radii $1$ and $\lambda$. The quotient
$\mathbb H^3/H_1$ is an open hyperbolic solid torus. As we adjoin the elements generated by $a_0$ to the
group, two thirds of the annulus becomes redundant,
a sector of $2\pi/3$ degrees becomes the fundamental region for $H_0$ as in Figure \ref{annulusoverthree}. The hyperbolic
quotient again becomes a solid torus, obtained from a Dehn twist.

\vspace{.05in}

\begin{figure}[!h] \bct
\includegraphics[width=.6\textwidth]{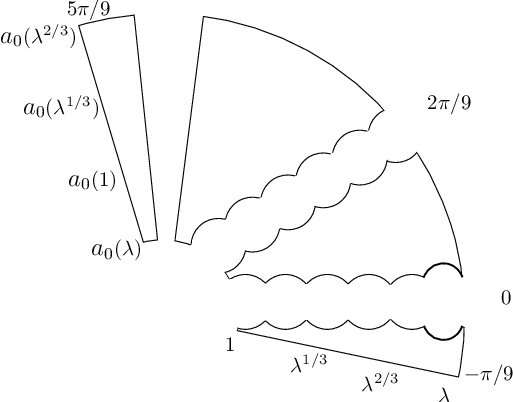}\\ \ect
  \caption{{\em {\small Fundamental region after $1/3$-complex twist.}}}
  \label{complextwistfundamentalregion}
\end{figure}

We have to normalize $G$ so that its fundamental region fits into the annulus piece. For this purpose,
 $G_2$ is joined into $G$ via conjugation $z\mapsto \exp(2\pi i/9)z$ by rotating $2\pi/9$ degrees rather than $2\pi/3$,
so that the identified circles stays inside the annular region between $-\pi/9$ and $5\pi/9$.
Besides, apply the first combination theorem to $G_1$ and $G_2$ taking the region $B_1$ as $|\arg z - 4\pi/9|<\varphi$
with $ \varphi < \pi/9 $,
and $B_2$ as before with its new lens angle $\varphi$. Now to combine the annular region with $G$, we take $B_1'$ as the annular region
$|-\pi/9 < \arg z < 5\pi/9|$ which is precisely invariant under $H_1=\<a\>$ in $H_0$. Take $B_2'$ to be the complementary region
$|5\pi/9 < \arg z < 2\pi -\pi/9|$ precisely invariant under $H_1$ in renormalized $G$.
Figure \ref{complextwistfundamentalregion} shows the resulting fundamental region.
Recall that $\mathbb H^3/H_0$ is a hyperbolic solid torus topologically obtained after applying
three Dehn twists to the solid torus $\mathbb H^3/H_1$. The ray $\{(z,t)| z=0 , t>0\}\subset \mathbb H^3$ projects onto the
central loop of the solid tori, where it is homotopic to the $(1,0)$ curve, the parallel on $\mathbb H^3/H_1$.
On the other hand it is homotopic to the $(1,3)$ torus knot on the boundary of $\mathbb H^3/H_0$. The second solid torus
is opened up along this homotopy, and glued back onto an opened up median of $\mathbb H^3/G$ in three dimensions.

\section{Handlebody Diagrams}\label{sechandlebody}


In this section we will draw  handlebody diagrams of the some
of the LCF $4$-manifolds constructed from the $3$ manifolds of the previous section via the application of the Theorem
\ref{braamthm}.
We will begin with $\Sigma_{1,2}$, the torus with two holes, then cross it with the interval $I=[0,1]$, and then glue the boundary cylinders with each other either trivially or with a flip. Then by gluing a solid torus to this (along the $p/q$ knot in its boundary)  to obtain the panelled web 3-manifold. We then cross this with $S^1$ and identify its boundary to obtain the panelled web 4-manifold.

\vspace{.05in}

\begin{figure}[!h]
\bct \includegraphics[width=0.5\textwidth]{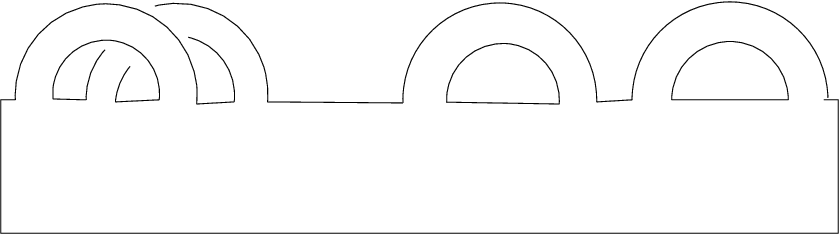}
 \caption{{\em {\small One handles of the torus with two punctures.} } } \label{2handlesoftorus}
\ect \end{figure}

\begin{figure}[!h]
\bct  \includegraphics[width=.3\textwidth]{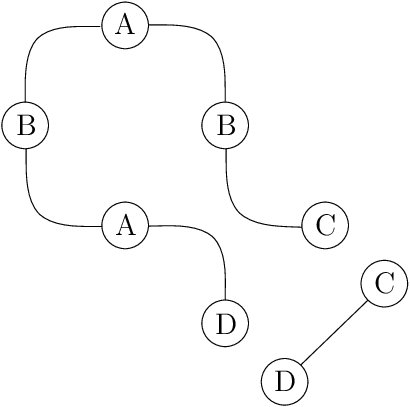}\\ \ect
  \caption{{\em {\small Heegard Diagram for $I\times
\Sigma_{1,2}$.} } } \label{heegard}
\end{figure}

Figure \ref{2handlesoftorus} is a handlebody picture of the twice punctured $2$-torus: It consists of  a $2$-disk (i.e. $0$-handle) with three $1$-handles attached to its boundary, and one $2$-handle (attached along the outer boundary of the figure).  Then Figure \ref{heegard} is just the thickening of this handlebody,
which is the {\em Heegard diagram} of $I\times \Sigma_{1,2}$.

\begin{figure}[!h]
\bct  \includegraphics[width=0.8\textwidth]{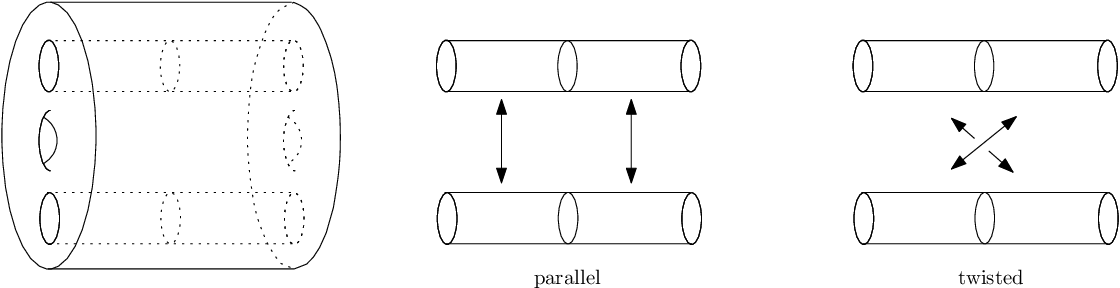}\\ \ect
  \caption{{\em {\small Identification of the boundary cylinders.} } } \label{cylinders}
\end{figure}

Now, we identify the two boundary cylinders in $I\times
\Sigma_{1,2}$ via the Second Combination Theorem of Maskit
\cite{maskitkleinian,maskitpwg}. We can do this in two different ways, either
trivially or with a twist. We will sketch the pictures of the
manifolds resulting from both ways of gluing.
This identification glues the neighborhoods of the middle
circles (called the {\em medians}\index{median} \cite{maskitpwg}) of the cylinders. As shown in Figure \ref{cylinders}.

\begin{figure}[!h]
\bct  \includegraphics[width=0.8\textwidth]{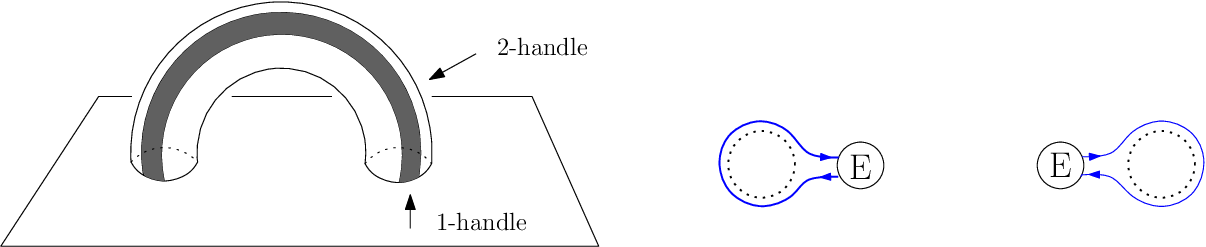}\\  \ect
  \caption{{\em {\small 2 and 3-dimensional round handles.} } } \label{23roundhandle}
\end{figure}

This operation of identifying the neighborhoods of the two
circles, is usually called the attaching a {\em round
$1$-handle}\index{round handle} operation. A round $1$-handle is a combination of a 1-handle and
a 2-handle as illustrated in Figure \ref{23roundhandle}.

\begin{figure}[!h] \bct
\includegraphics[width=0.8\textwidth]{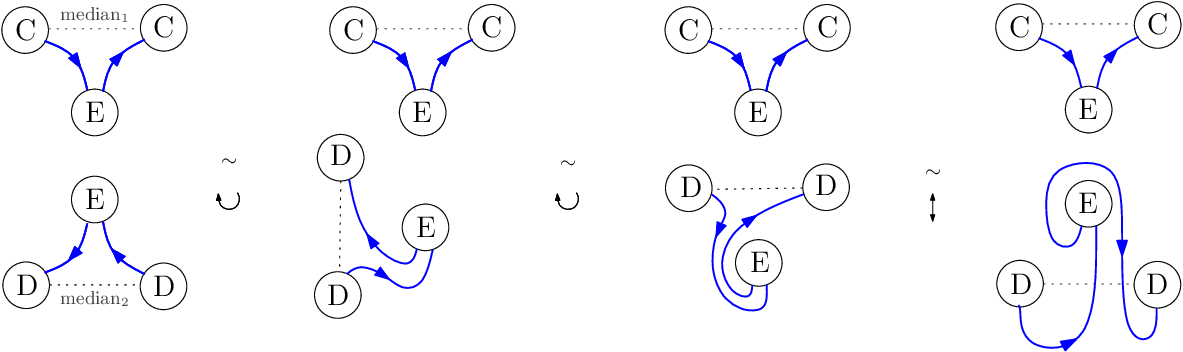}\\  \ect
  \caption{{\em {\small Obtaining the parallel round handle from the cross round handle.} } } \label{parallelcross}
\end{figure}

In the diagram of Figure \ref{parallelcross}, the median$_1$ and the median$_2$ are the cores of the $1$-handles C and D, respectively. This is because the median circles lie on the
cylinders, which make the holes on the 3-manifold, and we
formed these holes by the 1-handles C and D.

\vspace{.05in}

\begin{figure}[!h] \bct
\includegraphics[width=\textwidth]{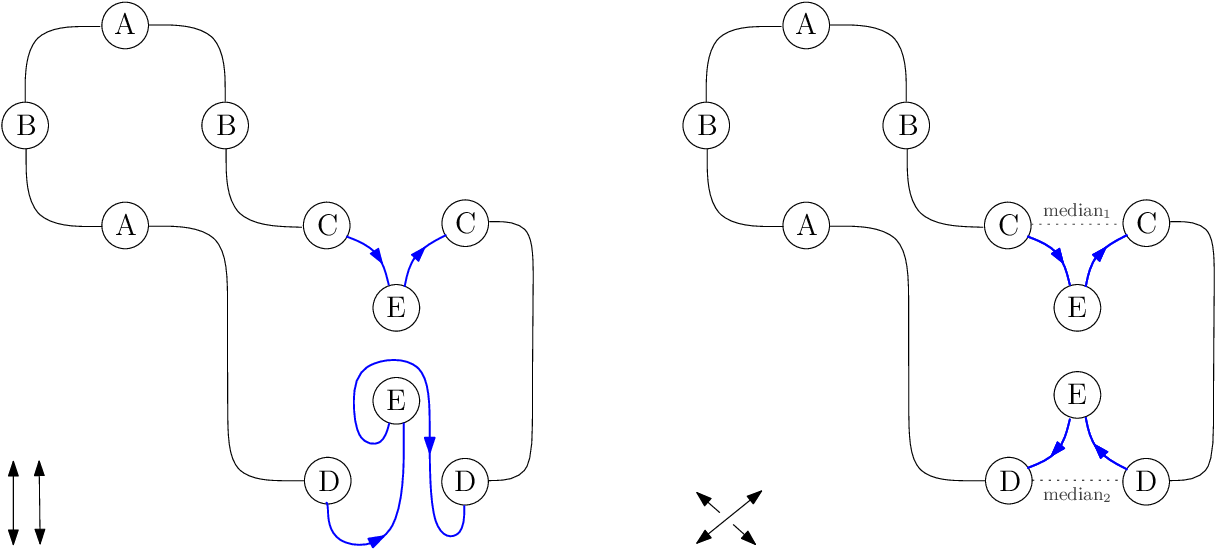}\\ \ect
  \caption{{\em {\small Two different ways of inserting the round handle.} } }\label{insertroundhandle}
\end{figure}

There are two different ways of gluing the neighborhoods of the meridians. Both ways are illustrated in Figure \ref{parallelcross}. In our figure we flipped the hole i.e. the 1-handle so that we can
obtain one identification from the other. We will call one
{\em cross identification} (the left picture), and the other  {\em parallel identification} (the right picture). In general the two different ways of attaching the round $1$-handles give non-diffeomorphic $3$-manifolds. (e.g. Figure \ref{insertroundhandle})

\begin{figure}[!h] \bct
\includegraphics[width=.25\textwidth]{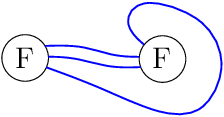}\\ \ect
  \caption{{\em {\small 1/3 torus knot on the 1-handle.} } } \label{1/3torusknot}
\end{figure}

\begin{figure}[!h] \bct
\includegraphics[width=.5\textwidth]{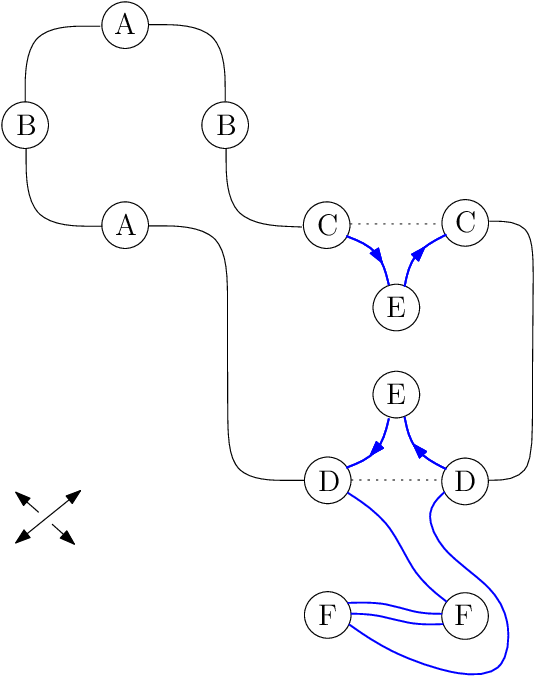}\\ \ect
  \caption{{\em {\small Maskit's 1/3 complex twist operation.} } } \label{1/3complextwist}
\end{figure}

The final operation to perform is to add a $p/q$ twist to this handlebody by gluing a solid torus to it.  This is done by identifying an annulus on its boundary with a neighborhood of a $p/q$ torus knot on the boundary of the solid torus, where
$p$ is the multiplicity of the meridian direction. Since the $p/q$ curve is isotopic to $1/q$ curve in the solid torus, it  suffices to take $p=1$. The solid torus here is viewed as a 1-handle, with  a $p/q$ torus knot lying on its boundary. In Figure \ref{1/3torusknot} we sketch the $1/3$ torus knot as an example.
This operation is similar to attaching a round handle operation (since we are identifying two circles), it is achieved with
a 1-handle and a 2-handle addition as in Figure \ref{1/3complextwist}. This finalizes
the picture of the Maskit's panelled web 3-manifold.

\vspace{.05in}
To pass  to  the 4-manifold, we cross this
3-manifold with a circle, and then shrink the boundary circles. Shrinking a circle is equivalent to identifying it to a point, which is achieved by attaching a 2-disk, we will call this {\em
capping the circle}\index{capping} operation.

\begin{figure}[!h] \bct
\includegraphics[width=.6\textwidth]{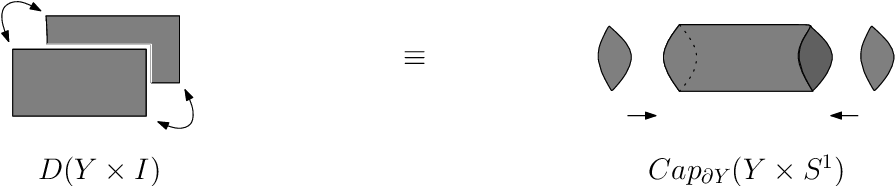}\\ \ect
  \caption{{\em {\small
  $D(Y\times I)=Cap_{\partial Y}(Y\times S^1)$ for the interval $Y$.} } }
  \label{thickening+double=crosscircle+capping}
\end{figure}

We begin by thickening the 3-manifold, i.e. crossing with an
interval. In particular, this amounts to thickening the pair of attaching $2$-disks of the three dimensional $1$-handle to $3$-balls (the attaching balls of the four dimensional $1$-handle).
The attaching circles
of the 2-handles inherit the blackboard framing from the 2-dimensional Heegard diagram. The blackboard framing can be computed as the {\em writhe}
\index{writhe} of the attaching knot of the 2-handle, i.e. the signed number of self crossings, which turns out to be 0 in our case. After thickening, we need to take the double of what we have. Thickening and taking the double is the same as
crossing with a circle and capping the boundary circles, as the lower
dimensional Picture \ref{thickening+double=crosscircle+capping}
illustrates.
Recall that the {\em double}\index{double} of a compact n-manifold
$X$ is defined to be
$$DX=\partial(I\times X)=X {\cup}_{
id_{\partial X} } \overline{X}.$$
where $\bar{X}$ is a copy of $X$ with the opposite orientation.
We denote the thickened 4-manifold by $X$, which is a 4-dimensional
handlebody without 3 or 4-handles. Then $DX$ automatically inherits
a handle decomposition: By turning the handle decomposition of $X$
upside down, we get the dual handle decomposition of $\overline{X}$,
which we attach on top of $X$ getting $DX=X \cup \mbox{dual
handles}$. Note that the duals of  0, 1 and 2-handles are 4, 3, and
2-handles, respectively. Since 3-handles are attached in a unique way, they don't need to be indicated in the picture.

\vspace{.05in}

\begin{figure}[!h] \bct
\includegraphics[width=.5\textwidth]{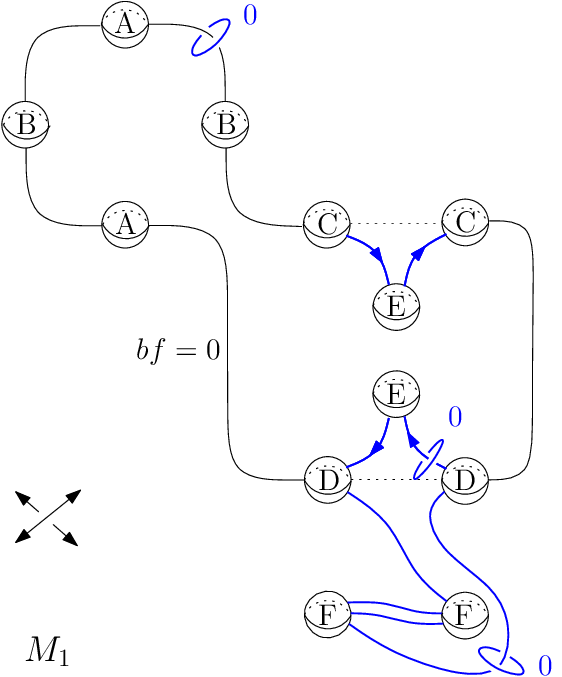}\\ \ect
  \caption{{\em {\small Thickening the Heegard Diagram and taking the double.} } } \label{thickanddouble}
\end{figure}

Hence to draw a handlebody picture of the double $DX$ from a given handlebody picture of $X$, it suffices to understand the position of the new (dual) 2-handles. They are attached by the $Id_{\partial X}$ map, along the cocores of the original 2-handles on the boundary.
So to get the double  we insert a 0-framed meridian
to each framed knot, as in the example in Figure \ref{thickanddouble}.
The 3 and 4-handles are attached afterwards uniquely to obtain the closed 4-manifold (they don't need to be drawn in the figure). We will denote this closed manifold by $M_1$\index{$M_1$}, it  corresponds the cross identification. We will denote the
manifold obtained from the parallel identification by
$M_2$\index{$M_2$}. Let us denote the corresponding manifolds (with boundary) before the doubling process, by $M_1'$ and $M_2'$ respectively, they only have  0,1 and 2 handles.



\vspace{.05in}

\begin{figure}[!h] \bct
\includegraphics[width=.4\textwidth]{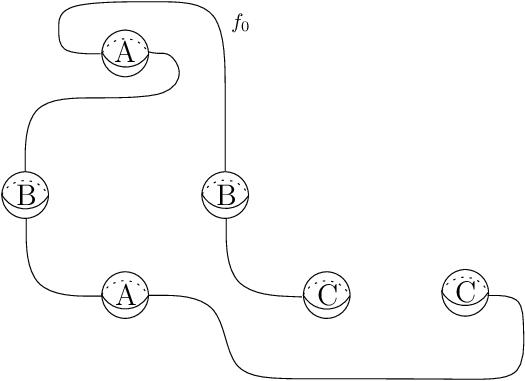}\\ \ect
  \caption{{\em {\small Twisted disk bundle over the punctured Klein bottle.}
  } }
  \label{thickeningtwistedklein1Ibundle}
\end{figure}

\begin{figure}[!h] \bct
\includegraphics[width=.8\textwidth]{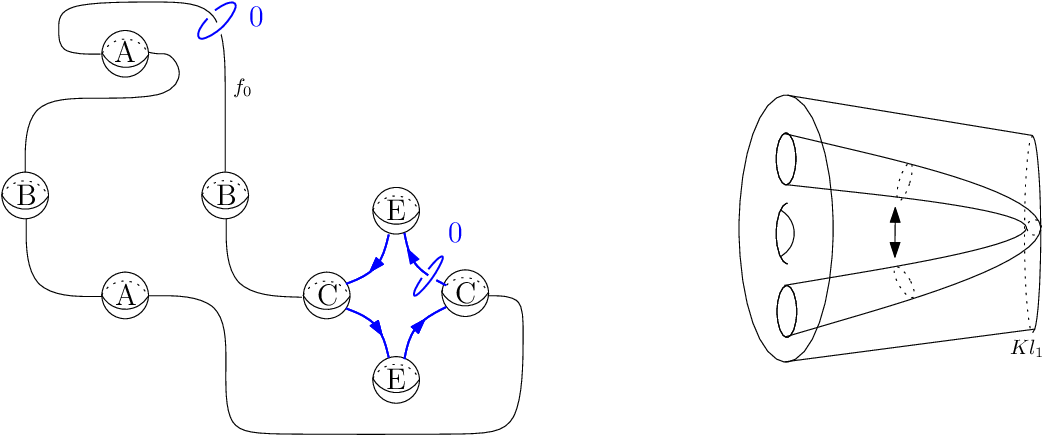}\\ \ect
  \caption{
{\em {\small Round handle and the double with the corresponding
3-manifold.}}}
  \label{klrouddouble}
\end{figure}

\begin{figure}[!h] \bct
\includegraphics[width=.6\textwidth]{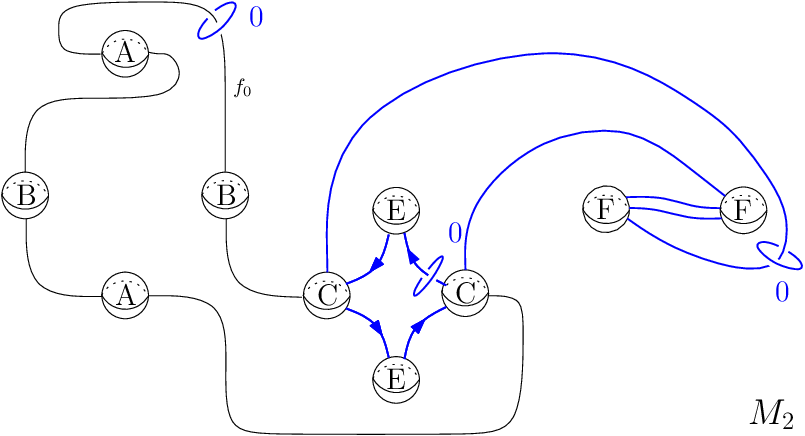}\\ \ect
  \caption{{\em {\small Maskit's 1/3 complex twist operation.}}}
  \label{1/3complextwist2}
\end{figure}

Now we treat  the twisted $I$-bundle case associated to the
surface $\Sigma_{1,2}$. Take a freely acting orientation reversing
involution $h:\Sigma_{1,2} \to \Sigma_{1,2}$, and
extend it to an orientation preserving homeomorphism
$$h':\Sigma_{1,2}\times I\to \Sigma_{1,2}\times I~~~\tn{by}~~~h'(x,t)=(h(x),1-t).$$
The resulting quotient $\Sigma_{1,2}\times I/h'$ is a twisted
$I$-bundle over a punctured Klein bottle $Kl_1$, which we denote by $Kl_1\wt{\times}I$. This could be thought as the quotient
$\Sigma_{1,2}\times I/ \sim$ as well, where $(x,1)\sim
(h(x),1)$.
Next we thicken and then double it. The thickening will result in
$Kl_1\wt{\times} I \times I \approx Kl_1\wt{\times}D^2$, a twisted
disk bundle over the punctured Klein bottle.
\begin{figure}[!h]
\bct \includegraphics[width=.4\textwidth]{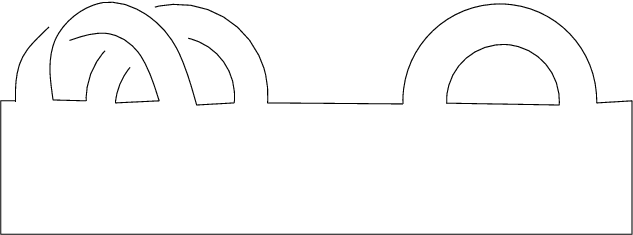}
 \caption{{\em {\small One handles of the punctured Klein bottle $Kl_1$.} } } \label{1handlesofpuncturedklein}
\ect \end{figure}
Figure \ref{1handlesofpuncturedklein} is the
handlebody of the punctured Klein bottle. Assuming that the framing
is the number $f_0$, the twisted disk bundle over the punctured Klein
bottle is sketched as in Figure
\ref{thickeningtwistedklein1Ibundle}.
Attaching the round handle $E$ and taking the double yields the
Figure \ref{klrouddouble}. Here, realize that there is a unique way
to attach the round handle according to Maskit's procedure. The
3-manifold is also drawn besides the 4-manifold picture.
Also, as before, we may twist by $1/3$ to obtain the Figure
\ref{1/3complextwist2}. We denote the resulting manifold by
$M_3$\index{$M_3$}, and the manifold with boundary before doubling
by $M_3'$.

\vspace{.05in}

\begin{figure}[!h] \bct
\includegraphics[width=\textwidth]{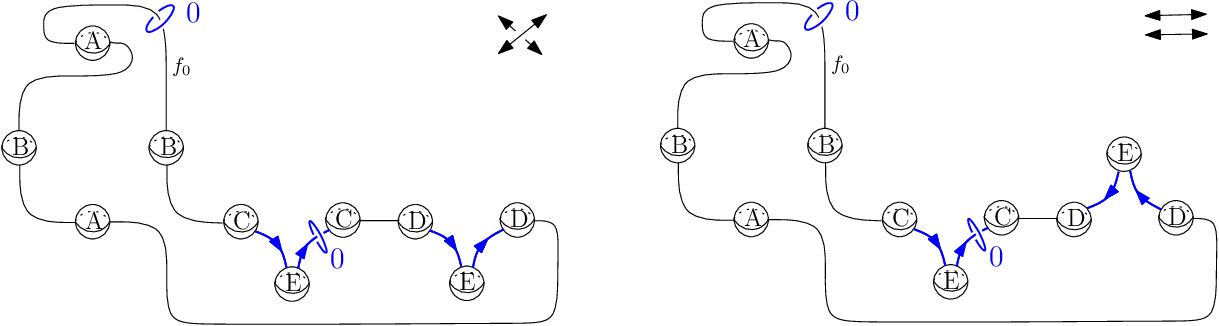}\\ \ect
  \caption{{\em {\small Cross and parallel identifications of the boundary cylinders of $Kl_2\wt{\times}D^2$.}}}
  \label{crossparalleltwicepuncturedklein}
\end{figure}

\begin{figure}[!h] \bct
\includegraphics[width=.6\textwidth]{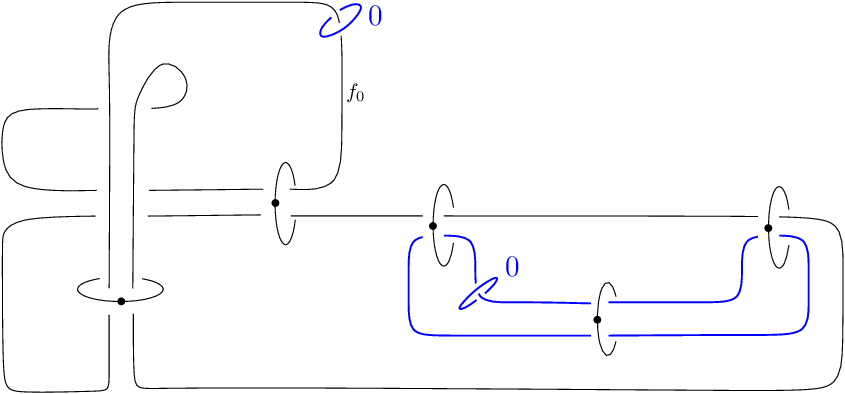}\\ \ect
  \caption{{\em {\small Dotted circle convention for the cross manifold of  Figure \ref{crossparalleltwicepuncturedklein}.}}}
  \label{dottedcircle}
\end{figure}

As a third example, we consider the twisted $I$-bundle over the twice punctured Klein bottle. We glue the boundary cylinders of
the twisted disk bundle over  $Kl_2$ in the cross and parallel
fashion to obtain the Figure \ref{crossparalleltwicepuncturedklein}.
After these operations, one may want to add the complex twists as well.  
To simplify the figures, one can use the dotted circle notation of
\cite{akbdot} to present our 4-manifolds. For example, Figure
\ref{dottedcircle} is the alternative handlebody picture of the
cross manifold just constructed. 

\vspace{.05in}

Here, we give a procedure of identifying the boundary cylinders of different manifolds. Note that whenever we draw two handlebody diagrams of 4-manifolds next to each other, it means that their handles
are attached on a common $S^3$ i.e. they have the same $0$-handle  $D^4$. So that they can be thought as two separate handlebodies  connected
by a 1-handle. Hence we only need to use the 2-handle
of the round handle to identify the two cylinders. This is how the
identification performed for the first pair of cylinders. For the rest of the identifications
the regular procedure applies, that is to build a tube (round handle) we need a 1-handle over which the 2-handle passes.

\vspace{.05in}

\begin{figure}[!h] \bct
\includegraphics[width=\textwidth]{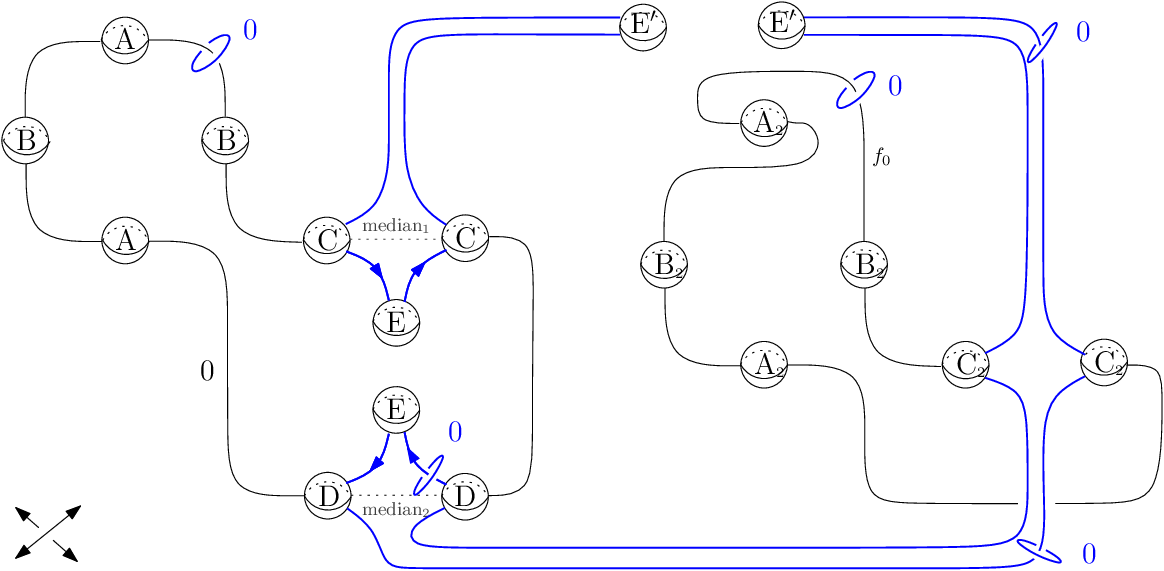}\\ \ect
  \caption{{\em {\small 4-manifold corresponding to the Maskit's  example.}}}
  \label{maskitmanifold}
\end{figure}

Finally, we draw the handlebody of the 4-manifold corresponding to an example of Maskit,
which he constructed from two different (trivial and twisted) types of
$I$-bundles associated to a torus with two holes namely $P_1$ and
$P_2$. He pairs the two ends of $P_2$ with a pair of cross ends of
$P_1$, the remaining cross ends of $P_1$ are identified with one
another. This  4-manifold is given by Figure
\ref{maskitmanifold}. 
Here, $E$ is the 1-handle of the round handle attaching a pair of
cross ends of the 4-manifold corresponding to $P_1$. Also  $C_2$ is identified to $D$ by
using only a 2-handle, and $E'$ is the 1-handle of the second round handle identifying
$C_2$ to $C$.

\section{Invariants}\label{secinvariants}

In this section we compute the topological invariants of the
manifolds constructed in the previous section. We first write down
the generators and relations of the fundamental groups.
We begin with the first set of construction (Figure
\ref{thickanddouble}). Each 1-handle is a generator of the
fundamental group, and each 2-handle provides a relation. We call
the generators $a,b,c,d,e,f$. We take the convention of left to
right and top to bottom to be the positive directions. Then, if we
begin from the
portion of the first 2-handle joining $D$ to $A$, going in the direction of $A$,
the first 2-handle provides the
relation
\beg{equation}\label{eq:torusrelation}a^{-1}b^{-1}abcd^{-1}=1.\end{equation}
If we begin with the 1-handle
$E$ of the round handle, its 2-handle gives the relation
\beg{equation}\label{eq:trround}   ede^{-1}c^{-1}=1.\end{equation}
Finally, the complex twist handle beginning with $F$ in the reverse
direction will provide 
\beg{equation}\label{eq:trcxtwist} f^{-3}d^{-1}=1.\end{equation}
If we abelianize this group, the first two relations yield the
relation $c=d$ and the third yields $c=f^{-3}$. Since $\< c,f ~|~
cf^3=1 \>=\<f^{-3},f\>=\<f\>$ the abelianization reduces the number
of generators by 2, hence $$H_1(M_1,\mathbb Z)=\mathbb Z^4.$$
Computing the second homology group needs more care. Since in the
doubling process we attach the upside down handles. Corresponding 
to each 1-handle, we have a 3-handle. So that the handles generate
the chain complex
$$0\to C_4 \to C_3 \to C_2 \to C_1 \to C_0 \to 0$$
$$0\to \mbb Z \to \mbb Z^6 \to \mbb Z^6 \to \mbb Z^6 \to \mbb Z \to 0.$$
This gives us the Euler characteristic $\chi(M_1)=1-6+6-6+1=-4$.
So  in terms of Betti numbers $-4=2b_0-2b_1+b_2$, implying
$b_2(M_1)=2$.
This is the free part. Next we compute the torsion piece.
By Poincar\'e duality $H_2(M_1,\mbb Z)\approx H^2(M_1,\mbb
Z)$, and since $H_1(M_1,\mbb Z) $ free the first term of the Universal Coefficients Theorem (e.g. \cite{hatcher}) is zero, we compute

$$0\to Ext(H_1(M_1,\mbb Z),\mbb Z) \to H^2(M_1,\mbb Z)\to Hom(H_2(M_1,\mbb Z),\mbb Z)\to 0$$
 $$H_2(M_1,\mbb
Z)=\mbb Z^2.$$
Similarly, we get $H_3\approx H^1 \approx H_1$   (by Poincare duality, and $H_0(M_1,\mbb Z)$ is free)
$$H_3(M_1,\mbb Z)=\mbb Z^4.$$

The alternative attachment of the round handle as $E$ in Figure
\ref{insertroundhandle}(a) gives the alternative for the second
relation (\ref{eq:trround}) \beg{equation}ed^{-1}e^{-1}c^{-1}=1\end{equation} which yields
$c=d^{-1}$ in the abelianization process, combining with the $c=d$
of (\ref{eq:torusrelation}) yields $c^2=1$. This implies that the
relation $d=f^{-3}$ of (\ref{eq:trround}) enforces $f^6=1$. So that
the first homology group becomes
$$H_1(M_2,\mathbb Z)=\< a,b,e,f ~|~ f^6=1 \> \approx\mathbb Z^3\oplus\mathbb Z_6.$$
The Euler characteristic $\chi(M_2)=-4$ since number of handles do not change, which implies $b_2(M_2)=0$. Also $Ext(\mbb Z^3\oplus \mbb Z_6,\mbb Z)=\mbb Z _6 $ becomes
the torsion part of $$H_2(M_2,\mbb Z)=\mbb Z_6.$$
Again by $H_3\approx H^1\approx Hom(H_1,\mbb Z)$ we have
$$H_3(M_2,\mbb Z)=\mbb Z^3.$$

Similarly, in the second set of constructions, in Figure
\ref{1/3complextwist2} we have the relations $$a^{-1}bab c=1 ~~,~~
ec^{-1}e^{-1}c=1 ~~,~~ f^{-3}c=1.$$ The first and third relation
imposes restrictions so that
$$H_1(M_3,\mbb Z)=\<a,b,c,e,f~|~c=b^{-2}=f^3\>=\<a,e,bf\>\approx\mbb Z^3$$
since $(bf)^{3}=b$, $(bf)^{-2}=f$ and $(bf)^{-6}=c$. The Euler
characteristic is $\chi(M_3)=1-5+6-5+1=-2$. So $b_2(M_3)=2$. $H_1$ and
$H_0$ has no torsion, hence
$$H_2(M_3,\mbb Z)=\mbb Z^2~~~\textnormal{and}~~~H_3(M_3,\mbb Z)=\mbb Z^3.$$

The signatures are $\sigma(M_{1,2,3})=0$ so that
$b_2^\pm(M_{1,3})=1,\, b_2^\pm(M_2)=0$ 
and the the intersection forms are \cite{braam}

$$Q_{M_{1,3}}=
\left[\beg{array}{cc}
0&1\\1&0\end{array}\right] :=H ~~~\textnormal{and}~~~Q_{M_2}= (0).$$ 
The invariants of the other two type of variations can be similarly
calculated.

\section{Sequences of Metrics}\label{secsequences}

Our goal in this section will be to combine our building blocks to construct
some interesting sequences of 4-manifolds admitting LCF metrics.
We begin by exploiting the first example described by Figure \ref{thickanddouble}.
There is no harm to replace the torus,
with any genus-g surface. We call the 4-manifolds arisen this way as $M^1_g$.
In this case the relation
$$a_1^{-1}b_1^{-1}a_1b_1\cdots a_g^{-1}b_g^{-1}a_gb_gcd^{-1}=1$$
replaces the relation (\ref{eq:torusrelation}); other relations (\ref{eq:trround}),(\ref{eq:trcxtwist}) remains. If we let
$g\longrightarrow\infty$, then we obtain
$$b_1(M_g^1)=2g+2\to\infty,$$ $$b_2(M_g^1)=2,$$ $$\chi(M_g^1)=-4g\to
-\infty.$$
Clearly $\sigma(M_g^1)=0$ and
$Q_{M_g^1}=H,$ both stay constant as we take the limit.

\vspace{.05in}

\begin{figure}[!h] \bct
\includegraphics[width=\textwidth]{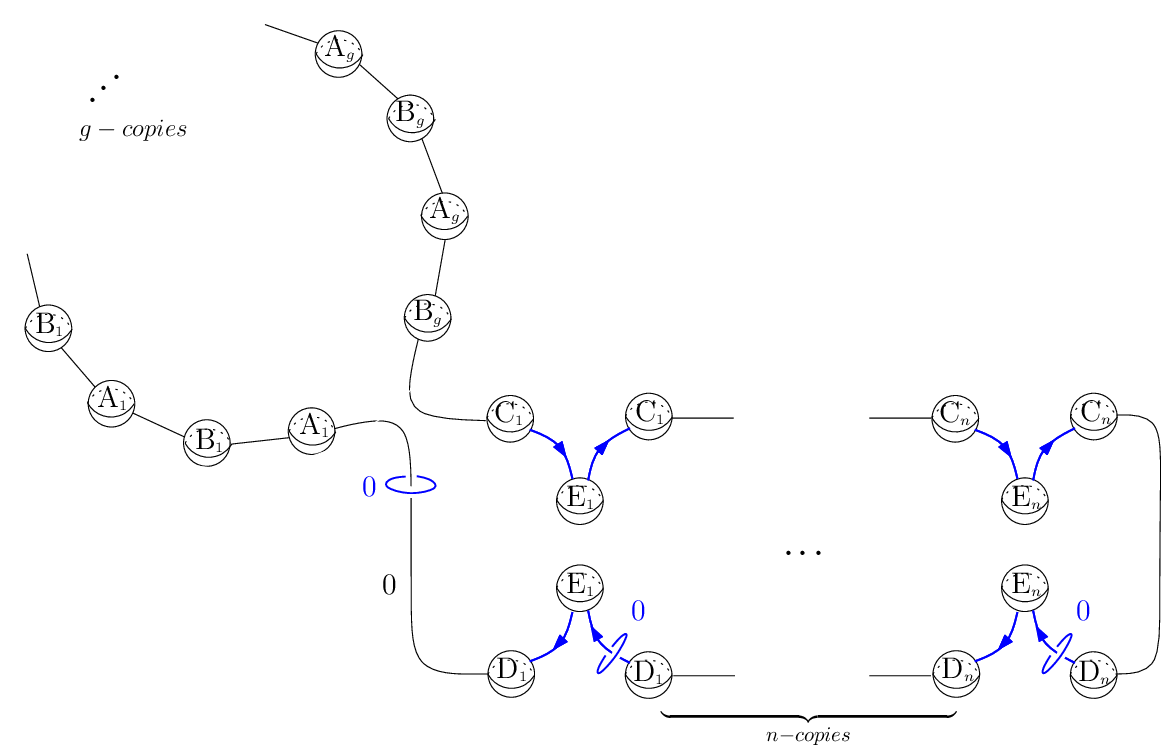}\\ \ect
  \caption{{\em {\small The LCF manifolds $M^2_{g,n}$.} } } \label{sequenceofmetrics12}
\end{figure}

Secondly, we may increase the number of $CDE$ components in (\ref{thickanddouble}) and
omit the complex twist 
handle $F$ for simplicity. We denote the resulting manifold $M^2_{g,n}$ 
(or sometimes $M^1_{g,n}$) where
$n$ stands for he number of $CDE$ components. See Figure \ref{sequenceofmetrics12}.
The orientations for $A$ handles are taken to be counterclockwise, and for $B$ handles to be clockwise.
The relations for 1-handles are
$$a_1^{-1}b_1^{-1}a_1b_1\cdots a_g^{-1}b_g^{-1}a_gb_gc_1\cdots c_n d_n^{-1}\cdots d_1^{-1}=1$$
$$e_i d_i e_i^{-1}c_i^{-1}=1 ~~\textnormal{for}~~ i=1\cdots n.$$
So that we obtain $$b_1(M^2_{g,n})=2g+2n \rightarrow\infty,$$  
$$b_2(M^2_{g,n})=2 $$ and $$\chi(M^2_{g,n})=4-4g-4n\to
-\infty $$ as $n\longrightarrow\infty$, and the intersection forms are
given by $Q_{M^2_{g,n}}=H$. \\

\begin{figure}[!h] \bct
\includegraphics[width=.4\textwidth]{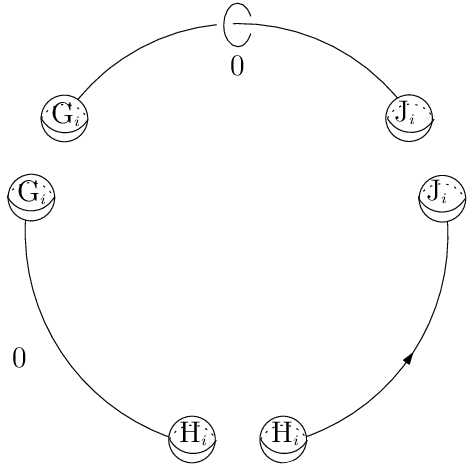}\\ \ect
  \caption{{\em {\small Doubling the $D^2 \times\Sigma_{0,3}$.} } } \label{doublingtheD2bundleover3puncturedS2}
\end{figure}

\begin{figure}[!h]
\bct
\includegraphics[width=\textwidth]{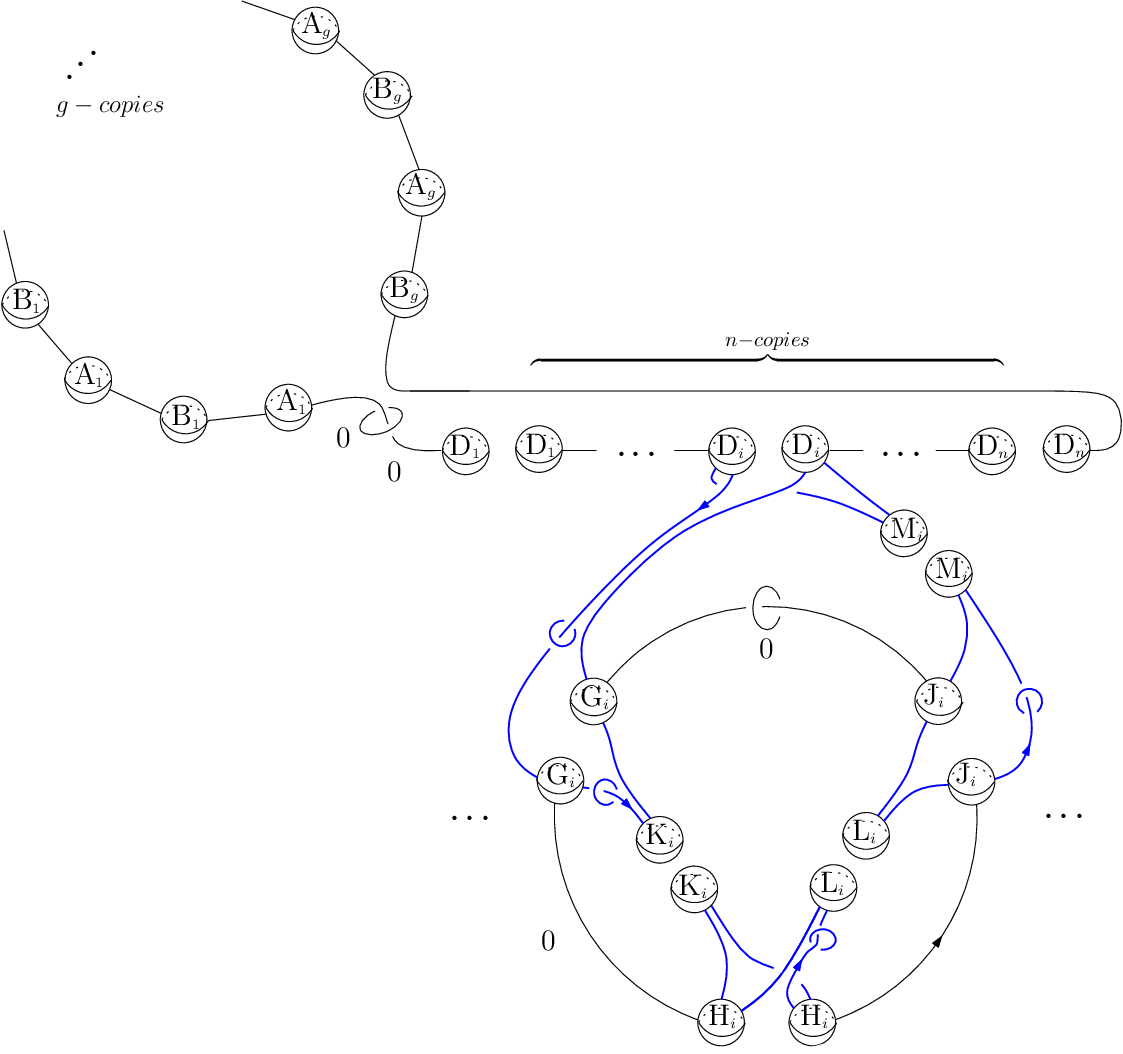}\\ \ect
  \caption{{\em {\small The LCF manifolds $M^3_{g,n}$.} } } \label{sequenceofmetrics3}
\end{figure}

In the third sequence, we will make use of another building block. This will be the
trivial I-bundle over a punctured annulus $\Sigma_{0,3}$. The corresponding 4-manifold
can be obtained by doubling the trivial disk bundle over $\Sigma_{0,3}$. Disk bundles
over $S^2$ are sketched as n-framed unknot. We only need to dig holes by attaching three 1-handles.
As a result the handlebody diagram is going to look as in Figure \ref{doublingtheD2bundleover3puncturedS2}.
We could have cancelled the 2-handles along with a 1-handle and this
makes it diffeomorphic to $S^1\times S^3 \sharp S^1\times S^3$. However
we cannot make any handle cancellation at this point as it will
destroy one of holes which we are using for attachment. Next we will
attach this piece through the $D_i$ handles. Since we are attaching
a different manifold, the round handle of the first identification
has no 1-handle, the rest of the round handles are as usual. We
attach it n-times and denote the resulting manifold by $M^3_{g,n}$.
See Figure \ref{sequenceofmetrics3}. The original 1-handle gives us
a similar relation
$$a_1^{-1}b_1^{-1}a_1b_1\cdots a_g^{-1}b_g^{-1}a_gb_g d_n^{-1}\cdots d_1^{-1}=1  ~~\Rightarrow~~ d_1\cdots d_n=1.$$
on the other hand each attached new piece provides the relations
$$g_i^{-1}d_i^{-1}=1 ~~\Rightarrow~~ d_i=g_i^{-1},$$
$$k_ih_ik_i^{-1}g_i=1 ~~\Rightarrow~~ h_i=g_i^{-1},$$
$$l_i^{-1} j_i l_i h_i = 1 ~~\Rightarrow~~ h_i=j_i^{-1},$$
$$m_i^{-1}d_i^{-1}m_ij_i^{-1}=1 ~~\Rightarrow~~ d_i=j_i^{-1},$$
$$g_ih_ij_i=1 ~~\Rightarrow~~ j_i=1,$$
\noindent where the right hand side of the arrows indicate the outcome in the abelianization process, so that we obtain $1=j_i=h_i=g_i=d_i$ and
the three free variables $k_i,l_i,m_i$ emerge from each attachment.
Counting these along with $a_i,b_i$ for $i=1\cdots g$ we have
$$b_1(M^3_{g,n})=2g+3n.$$
The Euler characteristic is computed at the chain level as
$$\chi(M^3_{g,n})=2-2(2g+7n)+(10n+2)=4-4g-4n.$$
From here we get 
$$b_2(M^3_{g,n})=2+2n.$$ So that $b_1,b_2\to\infty$ and $\chi\to-\infty$ as
$n\longrightarrow\infty$. The main difference of this sequence of
metrics from the previous ones is that $b_2$ gets arbitrarily large
rather than staying constant. If we let $g\longrightarrow\infty$
instead, then $b_1\to\infty$ , $\chi\to -\infty$ and $b_2=$constant,
a behaviour similar to the previous situations.


\vspace{.02in}

\begin{figure}[!h]
\bct
\includegraphics[width=.8\textwidth]{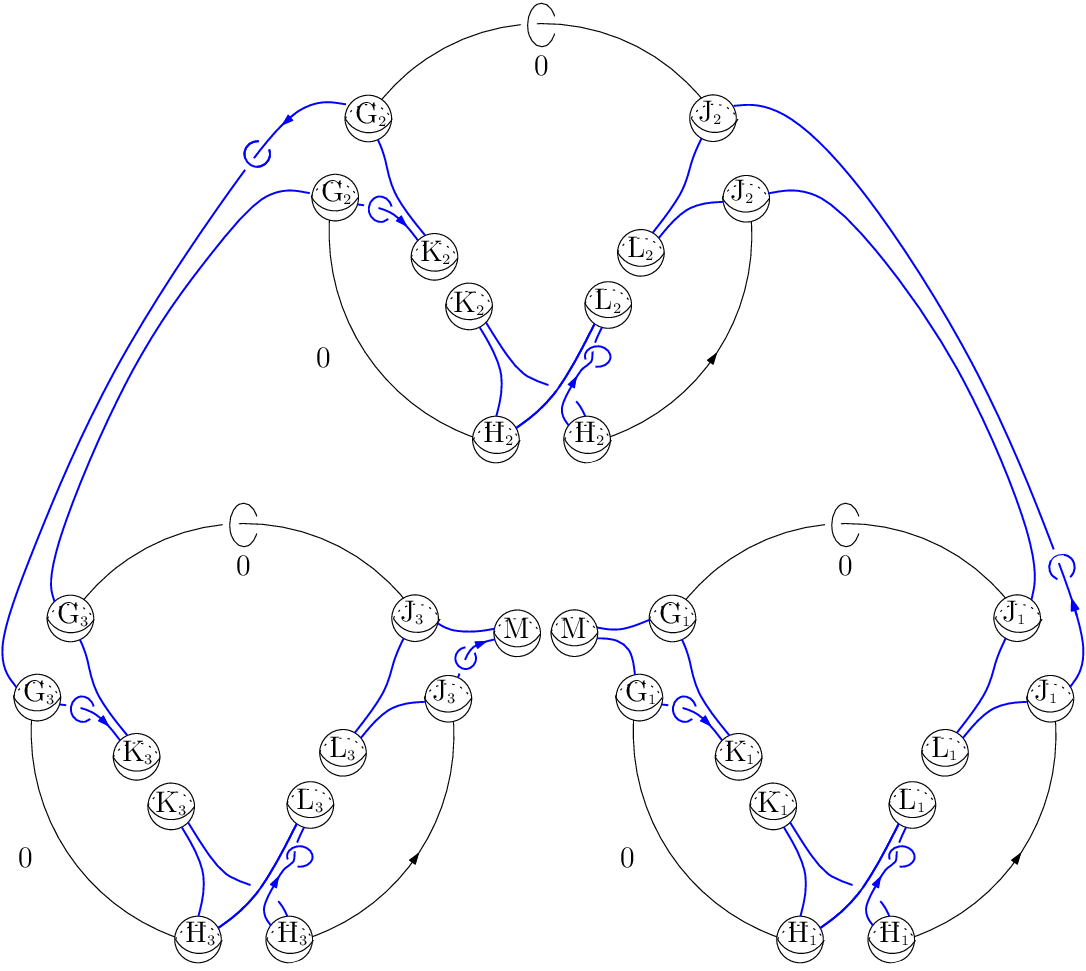}\\ \ect
  \caption{{\em {\small The LCF manifold $M^4_{3}$.} } } \label{sequenceofmetrics4}
\end{figure}

Our final sequence of panelled web manifolds is obtained by
attaching many copies of the new building block to each
other as a chain. One uses round handles without 1-handles to attach each copy,
and finally when closing up the line to a chain we use a complete
round handle. So that our chain contains only one complete round handle. 
Figure \ref{sequenceofmetrics4}
shows the case for $n=3$.
Again we have the relations
$$k_ih_ik_i^{-1}g_i=1 ~~\Rightarrow~~ h_i=g_i^{-1},$$
$$l_i^{-1} j_i l_i h_i = 1 ~~\Rightarrow~~ h_i=j_i^{-1},$$
$$g_ih_ij_i=1 ~~\Rightarrow~~ j_i=1.$$
The generators $g_i,h_i,j_i$ for the first homology are homologous to
each other and moreover are trivial. Only $k_i,l_i$ for $i=1\cdots
n$ and $m$ survive, so $$b_1(M^4_n)=2n+1.$$ The Euler characteristic
$$\chi(M^4_n)=2-2(5n+1)+8n=-2n,$$ and from these
$$b_2(M^4_n)=2n.$$
Again we have  $b_1,b_2\to\infty$ and $\chi\to-\infty$ as
$n\longrightarrow\infty$.

\section{Sign of the Scalar Curvature}\label{secsign}
In this section, we will verify the Theorem  
\ref{signscalarcurvaturenegative} on the sign of the scalar curvature. 
We will be using the results of LeBrun in \cite{cltopsd} in this section  unless otherwise stated. Main tool is the Weitzenb\" ock formula of 
\cite{bourguignon} involving the Weyl curvature. On a Riemannian manifold, the Hodge/modern Laplacian can be
expressed in terms of the connection/rough Laplacian as
$$(d+d^*)^2=\nabla^*\nabla-2W+{s\over 3}$$
where $\nabla$ is the Riemannian connection and $W$ is the Weyl
curvature tensor. First observation is that if there is a LCF metric of 
positive scalar curvature on a manifold, then the second Betti number 
$b_2=0$. Recall that  any de Rham cohomology class can be represented by a 
harmonic form uniquely on a closed manifold. One starts with an arbitrary 
harmonic 2-form and feeds it to the above formula. Then taking 
the inner product with the form and integrating over the manifold forces  the norm of the form to vanish. The zero scalar curvature case is more delicate. We will be using the following result, alternative exposition of which can be accessed through \cite{cloptimal} as well.
\begin{thm}[\cite{cltopsd}]
\label{SFASD4manifolds}
Let $(M,g)$ be a closed, scalar-flat anti-self-dual (SF-ASD) 4-manifold, 
then either \beg{itemize}
\item $b^+_2=0$, or
\item $b^+_2=1$ and $g$ is a scalar-flat K\" ahler metric; or else
\item $b^+_2=3$ and $g$ is a hyper-K\" ahler metric.
\end{itemize}
\end{thm}
\noindent The origin of the numbers $1$ and $3$ here is the possible number of the generating complex structures. Parallel self-dual 2-forms have constant length hence they correspond to compatible almost complex structures on a manifold and they are determined by their value at a point. Moreover each independent parallel form reduces the holonomy. 
If $b_2\neq 0$ for a SF-LCF manifold then since $\tau=0$ we are in the K\" ahler case. We are able to use the following result.
\begin{thm}[\cite{cltopsd}]
\label{kahlerSDspintypezero}
Let $(M,g)$ be a closed, self-dual, K\" ahler, spin 4-manifold of type zero, then $M$ is isometrically diffeomorphic to one of the following. \beg{itemize}
\item $K3$ surface with a Yau metric;
\item Flat 4-torus modulo a finite group; or
\item Flat 2-sphere bundle over a Riemann surface of genus $\geq 2$ with local product metric.
\end{itemize}
\end{thm}
\noindent The idea here is to use spin Weitzenb\" ock formula for nonzero signature to get a trivial canonical bundle, and in the zero signature case, the metric is LCF, 
and reducing the holonomy to a subgroup of $U(1)\times U(1)$ to get a Riemannian splitting.
Applying to our case, the first two cases are eliminated by the signature. 
panelled web 4-manifolds are not of the last two cases either. So that,  they are not of zero type either, in the $b_2\neq 0$ case. 

\vspace{.05in}

If one thinks in terms of metrics, one can verify this sign using the results of \cite{schoenyau88,nayatani97} even in the $b_2=0$ case.
Computation of the sign of the scalar curvature for our LCF  manifolds is related to the Hausdorff dimension of the Kleinian groups
used to uniformize the hyperbolic 3-manifold. A basic observation of  \cite{braam} is that the Kleinian
group $G$ of an hyperbolic 3-manifold acts on $S^4$ by the following orientation preserving conformal diffeomorphism:
$$i:\mbb H^3\times S^1 \to \mbb R^2\times (\mbb R^{2})^{*}\approx \mbb R^4-\mbb R^2\approx S^4-S^2$$ $$(x,y,t,\theta)\mapsto (x,y,t\cos\theta,t\sin\theta)$$
where $x,y\in\mbb R,t\in\mbb R^+$
are the coordinates of the hyperbolic space. The circle action in the domain
corresponds to the rotations of $\mbb R^2\times \mbb R^{2*}$ in the second component.
When we continuously extend this map to the boundary, we obtain the compactification map
$i:\overline{\mbb H}^3\times S^1 \to S^4.$
$\tn{PSL}(2,\mbb C)$ acts on $\overline{\mbb H}^3$ to result $\overline{M}^3$
as well as on $S^4$ on the right by conformal transformations, i.e. fractional linear
transformations
$$\mbb{HP}_1\times \tn{PSL}(2,\mbb C)\to\mbb{HP}_1$$
$$([x,y],\left[\begin{array}{cc} a&c\\ b&d\end{array}\right]) \mapsto [xa+yb,xc+yd].$$
\noindent 
The circle action is free in the interrior, its fixed point set is the boundary $S^2\times S^1$, which maps to the $S^2$ of the image $S^4$. $S^1\times \tn{PSL}(2,\mbb C)$ acts equivariantly with respect to $i$. If $\Lambda$ is the limit set of $G$, the limit set of the $G$-action on $S^4$ equals $i(\Lambda\times S^1)$, since
the circle action does not move the boundary $S^2$ this limit set is isomorphic to $\Lambda$. Summarizing $\Lambda \subset \mbb{CP}_1\subset\mbb{HP}_1$. Considering the inclusions $G\subset \tn{PSL}(2,\mbb C) \subset \tn{PGL}(2,\mbb H)$, we can state the result
of Schoen-Yau and the refinement of Nayatani which helps us to compute the sign
\begin{thm}[\cite{schoenyau88,nayatani97}]
\label{schoenyaunayatani}
Let $(X,[g])$ be a compact, LCF 4-manifold, which is uniformized by taking the quotient of $\Omega\subset S^4$
by the Kleinian group $G\subset \tn{PGL}(2,\mbb H)$ of conformal transformations of $\mbb{HP}_1$.
Let $g\in [g]$ be a metric (in the conformal class) of constant scalar curvature which exists by the solution of the
Yamabe Problem. Assume that  the limit set $\Lambda$ of $G$
is infinite, and the Hausdorff dimension $\dim (\Lambda ) > 0$.
Then the sign of the scalar curvature is equal to the sign of $1-\dim (\Lambda )$.
\end{thm}
\noindent We will see that the LCF manifolds
constructed in the previous sections are all of (strictly) negative scalar
curvature type. To be able to make use of Theorem \ref{schoenyaunayatani} we begin with a definition and cite some results in hyperbolic geometry.

\beg{defn}[\cite{camita}] A compact irreducible 3-manifold $M$ with
incompressible boundary is called a {\em generalized book of
$I$-bundles} if one may find a disjoint collection $A$ of essential
annuli in $M$ such that each component $R$ of the manifold obtained
by cutting $M$ along $A$ is either a solid torus, a thickened torus,
or homeomorphic to an I-bundle such that $\partial R \cap \partial
M$ is the associated $\partial I$-bundle.
\end{defn}
\noindent For a hyperbolic 3-manifold $(M,g)$, let $d(M,g)$ or $d(M)$  denote the Hausdorff dimension of the limit set of the discrete group which
acts on the hyperbolic space isometrically to give $(M,g)$ as the
quotient. By minimizing $d$ over all of the supporting hyperbolic structures, we obtain a topological invariant of $M$:
$$D(M):=\inf ~\{d(M,g)~|~g ~\tn{is a complete hyperbolic metric on}~ M \}.$$

\beg{thm}[\cite{bisjon}] Let $M$ be a compact, orientable,
hyperbolic 3-manifold. If $d(M)=1$ then $M$ is either a handlebody
or an $I$-bundle. (If $d(M)<1$ then $M$ is a handlebody or a
thickened torus.) \end{thm}

\beg{thm}[\cite{camita} Main Theorem II,Corollary 2.4] Let $M$ be a
compact, orientable, hyperbolizable 3-manifold which is not a
handlebody or a thickened
torus. 
Then $D(M)\geq 1$. \end{thm}
\noindent
If we combine these two theorems, we see that ~$d(M,g)>1$~ for our
hyperbolic metrics. So that ~$1-d < 0$, hence the scalar curvature
is strictly negative for our LCF 4-manifolds according to the Theorem
\ref{schoenyaunayatani}. We should keep in mind the equality $d(M,g)=\dim(\Lambda)$,
as explained prior to the theorem.

The theorem of Schoen-Yau and the refinement of Nayatani
is actually more general than what we have stated in Theorem \ref{schoenyaunayatani}, and it is valid for all dimensions $n\geq 3$. The group of conformal transformations of the n-sphere is the group of isometries
of the hyperbolic (n+1)-ball by the Liouville's theorem \cite{docarmo}. The isometry group of the
hyperbolic ball on the other hand is computed by considering it as the imaginary upper
unit sphere in the Minkowski space $\mbb R^{n+1,1}$. The transformations that preserve the
indefinite metric and the orientation happen to preserve the upper sheet of the hyperboloid \cite{docarmo,petersen}
so that $$\tn{Conf}(S^n)=\tn{Isom}(B^{n+1}_h)=SO^\uparrow (n+1,1).$$
Consequently, the uniformizing Kleinian group is a subgroup of this Lie group.
In the particular cases we have \cite{cloptimal} $$\tn{Conf}(\mbb{HP}_1)=\tn{PGL}(2,\mbb H)=SO^\uparrow (5,1)$$
$$\tn{Conf}(\mbb{CP}_1)=\tn{PSL}(2,\mbb C)=SO^\uparrow (3,1).$$
In the general case, the sign of the scalar curvature is equal to the sign of the quantity $$\frac{n}{2}-1-\dim(\Lambda).$$

\bigskip

{\small \beg{flushleft} \textsc{Mathematics Department, Michigan State University
}\\
\textit{E-mail address :} \texttt{\textbf{akbulut@math.msu.edu}}
\end{flushleft}  }

{\small \beg{flushleft} \textsc{Mathematics Department, University of Wisconsin at Madison}\\
\textit{E-mail address :} \texttt{\textbf{kalafat@math.wisc.edu}}
\end{flushleft}  }


\end{document}